\documentclass[a4paper]{amsart}

\usepackage[utf8]{inputenc}
\usepackage[english]{babel}
\usepackage{stmaryrd,mathtools,commath,mathrsfs}
\newcommand{\ontop}[2]{\stackrel{\mathclap{\normalfont\mbox{\footnotesize{#1}}}}{\ensuremath{#2}}} 
\usepackage{bbm} 
\usepackage{amssymb,amsthm}
\numberwithin{equation}{section}
\usepackage{float}
\usepackage{graphicx}
\usepackage{yfonts,verbatim,calrsfs,enumitem}
\usepackage{comment}
\usepackage{todonotes}

\usepackage{amsmath}
\usepackage{xcolor}
\definecolor{cornellred}{rgb}{0.7,0.11,0.11}

\usepackage{cite}
\usepackage[colorlinks=false, allcolors=cornellred, pdfpagelabels]{hyperref}

\allowdisplaybreaks

\theoremstyle{plain}
\newtheorem{theorem}{Theorem}[section]
\newtheorem{lemma}[theorem]{Lemma}
\newtheorem{proposition}[theorem]{Proposition}
\newtheorem{corollary}[theorem]{Corollary}

\theoremstyle{definition}
\newtheorem{definition}[theorem]{Definition}
\newtheorem*{definition*}{Definition}

\newtheorem{remark}[theorem]{Remark}

\newtheorem{notation}[theorem]{Notation}

\DeclareMathAlphabet{\mathcal}{OMS}{cmsy}{m}{n}

\newcommand{\Cs}{\ensuremath{\mathrm{C}^\ast}}
\newcommand{\1}{\ensuremath{\mathbf{1}}}
\newcommand{\id}{\ensuremath{\mathrm{id}}}
\newcommand{\Aut}{\ensuremath{\mathrm{Aut}}}
\newcommand{\Ad}{\ensuremath{\mathrm{Ad}}}
\newcommand{\N}{\ensuremath{\mathbb{N}}}

\newcommand{\R}{\ensuremath{\mathbb{R}}}

\newcommand{\M}{\ensuremath{\mathcal{M}}}
\newcommand{\U}{\ensuremath{\mathcal{U}}}
\newcommand{\K}{\ensuremath{\mathcal{K}}}
\newcommand{\C}{\ensuremath{\mathcal{C}}}

\newcommand{\e}{\ensuremath{\varepsilon}}
\newcommand{\f}{\ensuremath{\varphi}}
\newcommand{\bbu}{\ensuremath{\mathbbm{u}}}

\newcommand{\D}{\ensuremath{\mathcal{D}}}
\newcommand{\Q}{\ensuremath{\mathcal{Q}}}

\usepackage{tikz}
\usetikzlibrary{arrows,matrix,positioning,automata}
\usepackage{quiver}
\usepackage{pgfplots}

\setlist[enumerate,itemize]{noitemsep} 

\begin{document}

\title[Coronas and strongly self-absorbing \Cs-dynamics]{Corona algebras and strongly self-absorbing \Cs-dynamics}


\author{Xiuyuan Li}

\author{Matteo Pagliero}
\address{Department of Mathematics and Computer Science, University of Southern  
\linebreak
\phantom{----}Denmark, Campusvej 55, 5230 Odense, Denmark}
\email{pagliero@imada.sdu.dk}

\author{Gábor Szabó}
\address{Department of Mathematics, KU Leuven, Celestijnenlaan 200B, box 2400,\linebreak
\phantom{----}3001 Leuven, Belgium}
\email{gabor.szabo@kuleuven.be}
\thanks{}




\begin{abstract}
This article concerns the structure of \Cs-algebraic group actions induced on corona algebras from a given $\sigma$-unital \Cs-dynamical system over a locally compact group $G$.
We prove that such actions satisfy the so-called dynamical folding property, which generalizes a fundamental property observed for corona algebras in works of Manuilov--Thomsen and Phillips--Weaver.
We then focus on corona actions induced from $G$-\Cs-dynamics that are assumed to absorb a given strongly self-absorbing and unitarily regular $G$-action $\gamma$.
It is proved that these corona actions are $\gamma$-saturated, which is a stronger property than being separably $\gamma$-stable.
Conversely, if one assumes that the underlying \Cs-dynamics absorbs the trivial action on the compact operators, then $\gamma$-saturation of the corona action is equivalent to the original action being $\gamma$-absorbing.
These results are a dynamical version of recent work by Farah and the third-named author.
\end{abstract}

\maketitle

\setcounter{tocdepth}{1}
\tableofcontents

\section*{Introduction}
\renewcommand\thetheorem{\Alph{theorem}}

This article is a continued deep dive into the structural properties of corona algebras associated to \Cs-algebras with various absorption properties \cite{Farah23, FS24}.
Recall that the corona algebra of a \Cs-algebra $A$ is defined to be the quotient $\Q(A)=\M(A)/A$ of the multiplier algebra, which generalizes the construction of the Stone--{\v C}ech boundary for topological spaces.
Strongly self-absorbing \Cs-algebras and the algebras tensorially absorbing them \cite{TW07} play a crucial role in the structure and classification theory of \Cs-algebras.
In a similar vein, given a second-countable locally compact group $G$, the notion of a strongly self-absorbing $G$-action was introduced and studied for a similar purpose \cite{Sza18a, Sza18, Sza17}.
Even in the special case of a trivial $G$-action on a strongly self-absorbing \Cs-algebra $\D$, it can be instructive to study the notion of equivariant $\D$-stability for group actions on \Cs-algebras.
This and similar properties have already been exploited in various dynamical classification results as an intermediate technical step \cite{MS12, MS14, Sza18b, Sza19, Suz21, GS24, PS25}.

As was observed by Farah \cite{Farah23} and then subsequently refined together with the last-named author \cite{FS24} for \Cs-algebras, absorption of a strongly self-absorbing \Cs-algebra can often be detected at the level of the corona.
Given how most of the theory of strongly self-absorbing \Cs-algebras has some dynamical analog, this begs the question whether a similar behavior can be observed for group actions induced on corona algebras.
This article aims to present just such a dynamical analog of the aforementioned recent \Cs-algebraic results.
To formulate these, we begin by introducing a dynamical generalization of ``$D$-saturation''.
Recall that a unital \Cs-algebra $A$ is said to be $D$-saturated if for every separable \Cs-subalgebra $C\subseteq A$, there exists a unital embedding $D \hookrightarrow A$ whose image commutes with $C$.
When $A$ is equipped with an algebraic action $\alpha:G\curvearrowright A$, and $G$ acts on $D$ via a continuous action $\gamma$, we say that $\alpha$ is $\gamma$-saturated when the analogous property holds with respect to equivariant embeddings; we refer to Definition \ref{def:gamma-saturated} for the precise statement.
Our main result in its abridged form can be summarized as follows. 
For the full version, which holds in a more general setting, see Theorem~\ref{thm:saturation} and Corollary~\ref{cor:converse}.

\begin{theorem} \label{thm:A}
Let $\alpha:G\curvearrowright A$ be an action on a separable non-unital \Cs-algebra, $\gamma:G\curvearrowright\D$ a strongly self-absorbing, unitarily regular action, and equip $\Q(A)$ with the algebraic action induced by $\alpha$.\
If $\alpha$ is $\gamma$-absorbing, then $(\Q(A),\alpha)$ is $\gamma$-saturated.
Furthermore, if $\alpha$ tensorially absorbs the trivial action on the compact operators $\K$, then the converse holds as well.
\end{theorem}

The ``furthermore'' part in the above theorem is the comparably straightforward aspect of the main result deduced in the last section, since it can be obtained rather easily from the results and techniques that were introduced recently in \cite{FS24}.
The rest of the statement, however, is proved via a careful dynamical extension of the ideas behind the proof of \cite[Theorem A]{FS24}.
At least for actions of compact groups, this generalization is mostly technical in nature and can be achieved via the similar conceptual steps and extra bookkeeping.
This was accomplished in the first-named author’s master’s thesis \cite{Li-thesis}, which this article extends.
We would like to stress, however, that some new conceptual steps are needed in the proof of the statement for general \Cs-dynamics, in particular to cover arbitrary actions of groups $G$ that are non-compact.
This stems from the added difficulty that one may only identify $G$-actions modulo {\it cocycle} conjugacy for such groups in various contexts, which cannot be improved to genuine conjugacy in general.

In order to overcome this technical difficulty, we prove and make use of the so-called {\it dynamical folding property} for $G$-actions induced on corona algebras; see Definition~\ref{def:dynamical folding} for the details.
In somewhat oversimplified terms, a (not necessarily continuous) action $\alpha: G\curvearrowright A$ on a \Cs-algebra has the dynamical folding property if a given equivariant embedding $(D_0,\delta)\to (A,\alpha)$ of a separable system can be extended to an embedding of a surrounding system $(D,\delta)\supseteq(D_0,\delta)$ if and only if $D$ can be embedded in this way into the \Cs-algebra of continuous paths with values in $A$.
Without the dynamics (i.e.\ with $G=\{1\}$), this property has been observed for corona algebras by Manuilov--Thomsen \cite{MT04} and Phillips--Weaver \cite{PW07}; the term ``folding'' is directly inspired by a comment about this property made in the latter article.
This proved instrumental in Manuilov and Thomsen's proof that the $E$-theory of any \Cs-algebra pair $(A,B)$ can be obtained as the $KK$-group of the associated pair $(\mathcal{C}_0(\R)\otimes A, \Q(B\otimes\K))$.
Here we prove the following result, which is Theorem \ref{thm:dynamical folding} in the main body of the article.

\begin{theorem} \label{thm:B}
Let $\alpha:G\curvearrowright A$ be an action on a $\sigma$-unital \Cs-algebra.
Then $(\Q(A),\alpha)$ has the dynamical folding property.
\end{theorem}

While itself of independent interest, Theorem~\ref{thm:B} is used in the proof of Theorem~\ref{thm:A}.
As was mentioned above, when one aims to prove Theorem~\ref{thm:A} for a group $G$ that is not necessarily compact, one can generalize most of the steps from the $\Cs$-algebraic proof, but ultimately one faces a cocycle obstruction in carrying out the last part of the proof.
Due to the specific context here, the theory of strongly self-absorbing actions can be used to see that this cocycle obstruction can be at least asymptotically trivialized with a continuous path of unitaries.
As a result of the dynamical folding property, this allows us to genuinely trivialize the given cocycle obstruction in the corona and hence we can obtain Theorem~\ref{thm:A}.

Apart from the application to our main result, the dynamical folding property has another immediate consequence.
In the context of Theorem~\ref{thm:A},  we show that all equivariant unital embeddings $(\D,\gamma)\to (\Q(A),\alpha)$ are mutually $G$-unitarily equivalent, i.e., conjugates of each other via unitaries in $\Q(A)$ that are fixed by $\alpha$; see Corollary \ref{cor:G-unitarily-equivalence}.

\subsection*{Acknowledgments.}
M.P.\ was supported by PhD-grant 1131625N funded by the Research Foundation Flanders (FWO), and by the Carlsberg Foundation grant CF24-1144.
G.S.\ was supported by research project G085020N funded by the Research Foundation Flanders (FWO), and the European Research Council under the European Union's Horizon Europe research and innovation programme (ERC grant AMEN--101124789).
Both M.P.\ and G.S.\ would like to thank Ilijas Farah for several inspiring discussions related to this article.

For the purpose of open access, the authors have applied a CC BY public copyright license to any author accepted manuscript version arising from this submission.


\section{Preliminaries on strongly self-absorbing actions}
\renewcommand\thetheorem{\arabic{section}.\arabic{theorem}}

In this section, we recall the necessary background material on strongly self-absorbing \Cs-dynamical systems.

\begin{notation}
Throughout the article $G$ denotes a second-countable, locally compact group.\
We will often (but not always) assume that an action $\alpha: G\curvearrowright A$ of $G$ on a \Cs-algebra is point-norm continuous.\
If it is not clear from the context, we will explicitly say that $\alpha$ is \textit{continuous} to avoid confusion, or call it an \textit{algebraic} action when we want to stress that continuity is not assumed.
Given an action $\alpha:G\curvearrowright A$ (continuous or not), we denote by $A^\alpha\subseteq A$ the  fixed point algebra.
Given an inclusion $A\subseteq B$ of \Cs-algebras, we will use the notation
\[
B\cap A' = \{ x\in B \mid xa=ax \text{ for all } a\in A\}
\]
and
\[
B\cap A^\perp = \{ x\in B \mid xa=ax=0 \text{ for all } a\in A\}.
\]
\end{notation}

The work in this article appeals to the categorical framework of \cite{Sza21} in some instances.\
In this context, the appropriate notion of isomorphism between actions is defined as follows.

\begin{definition}
Let $\alpha:G\curvearrowright A$ and $\beta:G\curvearrowright B$ be actions on \Cs-algebras.\
A \textit{cocycle conjugacy} from $\alpha$ to $\beta$ is a pair $(\f,\bbu)$ consisting of an isomorphism $\f:A\to B$ and a strictly continuous map $\bbu:G\to\U(\M(B))$ such that $\bbu_{gh}=\bbu_g\beta_g(\bbu_h)$, i.e., $\bbu$ is a $\beta$-cocycle, and $\f \circ \alpha_g = \Ad(\bbu_g) \circ \beta_g \circ \f$ for all $g,h\in G$.

If $\bbu_g=\1$ for all $g\in G$,  $\f$ is an equivariant isomorphism, or in other words, a \textit{conjugacy}.
Whenever there exists a (cocycle) conjugacy between $\alpha$ and $\beta$, one says that they are (\textit{cocycle}) \textit{conjugate}.
\end{definition}

Strongly self-absorbing actions are the dynamical analog of strongly self-absorbing \Cs-algebras (see \cite{TW07}).
We recall their definition from \cite[Definition 5.3]{Sza21}, which (as per a comment therein) replaces \cite[Definition 3.1]{Sza18a}.

\begin{definition}[{see \cite[Definition 5.3]{Sza21}}]\label{def:ssa}
Let $\D$ be a separable, unital \Cs-algebra equipped with an action $\gamma:G\curvearrowright\D$.\
One says that $\gamma$ is \textit{strongly self-absorbing} if there exist a cocycle conjugacy $(\f,\bbu)$ from $\gamma$ to $\gamma\otimes\gamma$ and a sequence of unitaries $u_n\in\U(\D\otimes\D)$ such that
\begin{align*}
& \lim_{n\to\infty}\| \f(d) - \Ad(u_n)(d\otimes\1_{\D}) \|=0,
& \lim_{n\to\infty}\max_{g\in K}\| \bbu_g - u_n(\gamma\otimes\gamma)_g(u_n)^* \|=0
\end{align*}
for all $d\in\D$ and every compact subset $K\subseteq G$.
\end{definition}

\begin{notation}
Suppose we are given a strongly self-absorbing action $\gamma$ as above.
A continuous action $\alpha: G\curvearrowright A$ on any \Cs-algebra is called \textit{$\gamma$-absorbing}, or \textit{$\gamma$-stable}, if $\alpha$ is cocycle conjugate to $\alpha\otimes\gamma$.
In the special case where $\gamma=\id_{\D}$, the trivial action on a strongly self-absorbing \Cs-algebra $\D$, we say that $\alpha$ is equivariantly $\D$-stable.
\end{notation}

An action is strongly self-absorbing precisely when it is \textit{semi-strongly self-absorbing} in the sense of \cite[Definition 4.1]{Sza18a}.
One can conclude the following as a consequence of \cite[Remark 1.11, Theorem 4.7]{Sza18a}.

\begin{theorem} \label{thm:ssa-action}
Let $\alpha:G\curvearrowright A$ be an action on a separable \Cs-algebra, and $\gamma:G\curvearrowright\D$ a strongly self-absorbing action on a separable unital \Cs-algebra.\
The following are equivalent:
\begin{enumerate}[leftmargin=*,label=\textup{(\roman*)}]
\item $\alpha$ is $\gamma$-absorbing,
\item there exists a sequence of unital $\ast$-homomorphisms $\f_n:\D\to\M(A)$ such that
	\begin{itemize}
	\item $\|[\f_n(d),a]\|\xrightarrow{n\to\infty}0$ for all $a\in A$ and $d\in\D$, and
	\item $\alpha_g(\f_n(d))-\f_n(\gamma_g(d)) \xrightarrow{n\to\infty}0$ uniformly on compact subsets of $G$ in the strict topology for all $d\in\D$.
	\end{itemize}
\item there exists an equivariant unital $\ast$-homomorphisms
\[
\f:\D\to (\M(A)_\infty\cap A')/(\M(A)_\infty\cap A^\perp).
\]
\end{enumerate}
\end{theorem}

Next we recall the concept of unitary regularity.
It may be viewed as a dynamical analog of the assumption on a (unital) \Cs-algebra $A$ that the quotient group $\U(A)/\U_0(A)$ is abelian\footnote{Here, $\U_0(A)$ denotes the connected component of the unit.}.
This is automatic for $K_1$-injective \Cs-algebras, for instance, which includes all strongly self-absorbing \Cs-algebras or \Cs-algebras absorbing them by a result of Winter \cite{Winter11}.

\begin{definition}[{see \cite[Definition 2.17]{Sza18}}]
Let $\gamma:G\curvearrowright D$ be an action on a unital \Cs-algebra.\
One says that $\gamma$ is \textit{unitarily regular} if for every compact subset $K\subseteq G$ and $\e>0$, there exists $\delta>0$ such that, whenever $u_1,u_2\in\U(D)$ satisfy $\max_{j=1,2} \max_{g\in K}\|\gamma_g(u_j)-u_j\|\leq\delta$, there exists a norm-continuous path of unitaries $(w_t)_{t\in[0,1]}\subseteq\U(D)$ such that
\begin{equation*}
 \max_{g\in K} \max_{t\in[0,1]}\|\gamma_g(w_t)-w_t\|\leq\e, \quad w_0=\1,\quad  w_1=u_1u_2u_1^*u_2^*.
\end{equation*}
\end{definition}

We stress that unitary regularity is a relatively mild condition; in fact it is currently unknown if strongly self-absorbing actions can fail this property.\
Note that by \cite[Proposition 2.19]{Sza18}, all equivariantly $\mathcal{Z}$-absorbing actions on unital \Cs-algebras are unitarily regular, which covers many examples.

It turns out that for strongly self-absorbing actions, unitary regularity can be characterised by a tensorial absorption property using the following auxiliary object.

\begin{definition}
Given a unital \Cs-algebra $D$, let us define
\begin{equation*}
D^{(2)} = \{f\in\C([0,1],D \otimes_{\max} D) \mid f(0)\in D\otimes\1, f(1)\in \1\otimes D\}.
\end{equation*}
Consider c.p.c.\ order zero maps $\eta_i: D\to D^{(2)}$ for $i=0,1$ given by
\begin{equation*}
\eta_i(d)(t) =
\begin{cases}
(1-t) (d\otimes\1) & i=0, \\
t (\1\otimes d) & i=1
\end{cases}
\end{equation*}
for all $t\in[0,1]$ and $d\in D$.
\end{definition}

\begin{remark}\label{rem:universal property}
Let $D$ be a unital \Cs-algebra.\
As a consequence of \cite[Lemma 5.2]{HWZ15} and \cite{WZ09} (see also \cite[Lemma 6.6]{HSWW17}), $D^{(2)}$ and the c.p.c.\ order zero maps $\eta_0$ and $\eta_1$ satisfy the following universal property.\
For every unital \Cs-algebra $B$ and c.p.c.\ order zero maps $\mu_i:D\to B$ for $i=0,1$ with commuting ranges and such that $\mu_0(\1)+\mu_1(\1)=\1$, there exists a unique unital $\ast$-homomorphism $\f:D^{(2)} \to B$ such that $\f \circ \eta_i = \mu_i$ for $i=0,1$.

As a result, if $\gamma:G\curvearrowright D$ is an action, it follows that there exists a well-defined continuous action $\gamma^{(2)}:G\curvearrowright D^{(2)}$ induced by $\gamma$ via the identity $\gamma^{(2)}_g \circ \eta_i = \eta_i \circ \gamma_g$ for all $g\in G$ and $i=0,1$.\
Note that $\gamma^{(2)}$ is the restriction of the action on $\C([0,1],D \otimes_{\max} D)$ that acts fiberwise by $\gamma\otimes\gamma$.
Moreover, it follows that the $G$-\Cs-algebra $(D^{(2)},\gamma^{(2)})$ has the following universal property.\
Whenever $(B,\beta)$ is a unital $G$-\Cs-algebra, and $\mu_i:D\to B$ for $i=0,1$ are equivariant c.p.c.\ order zero maps with commuting ranges and such that $\mu_0(\1)+\mu_1(\1)=\1$, then there exists a unique unital equivariant $\ast$-homomorphism $\f:(D^{(2)},\gamma^{(2)})\to(B,\beta)$ such that $\f \circ \eta_i = \mu_i$ for $i=0,1$.
\end{remark}

We are ready to recall from \cite{Sza18} the characterisation of strongly self-absorbing actions that are additionally unitarily regular, in terms of a property of the \Cs-dynamical system $(\D^{(2)},\gamma^{(2)})$.

\begin{theorem}[{see \cite[Theorem 5.9]{Sza18}}]\label{thm:unitarily-regular}
Let $\gamma:G\curvearrowright\D$ be a strongly self-absorbing action on a separable unital \Cs-algebra.\
Then $\gamma$ is unitarily regular if and only if $\gamma^{(2)}$ is $\gamma$-absorbing.
\end{theorem}

The path algebra associated to a \Cs-algebra is somewhat reminiscent of the (well-known) sequence algebra, and one may think of it as a continuous version of the latter.\
This object will play a central role throughout the article.

\begin{definition}
Let $A$ be a \Cs-algebra.\
Denote by $\C_b([0,\infty),A)$ the \Cs-algebra of bounded continuous functions from $[0,\infty)$ to $A$ and by $\C_0([0,\infty),A)$ its ideal of functions vanishing at infinity.
We denote by
\begin{equation*}
A_{\mathfrak{c}} = \C_b([0,\infty),A)/\C_0([0,\infty),A)
\end{equation*}
the \textit{path algebra} of $A$.\footnote{Note that the path algebra does not have a commonly agreed upon notation in the literature.}
Note that $A$ can be naturally identified with the \Cs-subalgebra of constant paths in $A_{\mathfrak{c}}$.\
When we need to make this embedding explicit, we use the notation $\iota_{\mathfrak{c}}$.
We moreover call by $\pi_{\mathfrak{c}}$ the quotient map $\C_b([0,\infty),A) \to A_{\mathfrak{c}}$.
\end{definition}

\begin{remark}\label{rem:inner-aut1}
Let $\sigma\in\Aut(A)$ be an automorphism of a \Cs-algebra.\
Then, by sending each $f\in\C_{b}([0,\infty),A)$ to the function $\sigma \circ f\in\C_{b}([0,\infty),A)$, one obtains an automorphism $\sigma_b$ of $\C_{b}([0,\infty),A)$.
This induces an automorphism $\sigma_{\mathfrak{c}}$ of $A_{\mathfrak{c}}$ because $\sigma \circ f$ vanishes at infinity when $f$ does.
The resulting assignments $\Aut(A)\to\Aut\big( \C_{b}([0,\infty),A) \big)$ and $\Aut(A)\to\Aut(A_{\mathfrak{c}})$ are clearly both multiplicative.
Moreover, if $C\subseteq A_{\mathfrak{c}}$ is a $\sigma_{\mathfrak{c}}$-invariant \Cs-subalgebra, then $\sigma_{\mathfrak{c}}$ restricts to an automorphism of $A_{\mathfrak{c}}\cap C'$ as well.
Note that if $A$ is unital, any inner automorphism of $A$ induces the identity map on $A_{\mathfrak{c}}\cap A'$.
\end{remark}

\begin{definition}
Let $\alpha:G\curvearrowright A$ be an algebraic action on a \Cs-algebra.\
Then we denote by $A_{\alpha}$ the subset of $A$ containing all elements $a$ on which $\alpha$ is continuous, i.e.,
\begin{equation*}
A_{\alpha} = \{ a\in A \mid [g\mapsto\alpha_g(a)] \text{ is continuous}\}.
\end{equation*}
This is a \Cs-subalgebra of $A$, which we refer to as the \textit{($\alpha$-)continuous part of $A$}.
\end{definition}

\begin{definition}
Let $\alpha:G\curvearrowright A$ be an action on a \Cs-algebra.
In light of the above, $\alpha$ induces algebraic $G$-actions $\alpha_b$ on $\C_{b}([0,\infty),A)$ and $\alpha_{\mathfrak{c}}$ on $A_{\mathfrak{c}}$.
In general these actions may fail to be continuous.
However, they restrict to continuous actions on their respective continuous parts, which we denote by
\begin{equation*}
	\C_{b,\alpha}([0,\infty),A) := (\C_b([0,\infty),A))_{\alpha_b} \quad\text{and}\quad A_{\mathfrak{c},\alpha} := (A_{\mathfrak{c}})_{\alpha_{\mathfrak{c}}}.
\end{equation*}
Note that as a \Cs-subalgebra of $A_{\mathfrak{c}}$, $A_{\mathfrak{c},\alpha}$ agrees with $\C_{b,\alpha}([0,\infty),A)/\C_0([0,\infty),A)$ thanks to \cite[Theorem 2]{Bro00}.
\end{definition}

\begin{proposition} \label{prop:embeddingD2c}
Let $\alpha:G\curvearrowright A$ be an action on a separable, unital \Cs-algebra, and $\gamma:G\curvearrowright \D$ a strongly self-absorbing, unitarily regular action on a separable, unital \Cs-algebra.\
The following are equivalent:
\begin{enumerate}[label=\textup{(\roman*)}]
\item There exists an equivariant unital embedding $(\D,\gamma)\hookrightarrow(A_{\mathfrak{c},\alpha}\cap A',\alpha_{\mathfrak{c}})$. \label{prop:embeddingD2c:1}
\item $\alpha$ is $\gamma$-absorbing. \label{prop:embeddingD2c:2}
\item There exists a cocycle conjugacy
\[
(\f,\bbu): (A,\alpha)\to(A\otimes\D,\alpha\otimes\gamma)
\]
and a continuous path of unitaries $w_t: [0,\infty)\to \U(A\otimes\D)$ with $w_0=\1$ such that
\[
0=\lim_{t\to\infty} \|\f(x)-w_t(x\otimes\1)w_t^*\|+\max_{g\in K} \|\bbu_g-w_t(\alpha_g\otimes\gamma_g)(w_t)^*\|
\]
for all $x\in A$ and every compact set $K\subseteq G$. \label{prop:embeddingD2c:3}
\end{enumerate}
\end{proposition}
\begin{proof}
\ref{prop:embeddingD2c:1}$\Rightarrow$\ref{prop:embeddingD2c:2}:
As $A_\mathfrak{c}\cap A'$ maps equivariantly into $A_\infty\cap A'$ by restricting functions to $\mathbb N\subset[0,\infty)$, this is a consequence of \cite[Theorem 4.7]{Sza18a}.

\ref{prop:embeddingD2c:2}$\Rightarrow$\ref{prop:embeddingD2c:3}:
In order to prove this implication, one follows the argument towards \cite[Theorem 3.2]{Sza17} word for word.
Said theorem was stated with a formally weaker conclusion, but one obtains the statement of \ref{prop:embeddingD2c:3} via the same argument.
In its proof, one verifies that the equivariant embedding
\[
\1_{\D}\otimes\id_A : (A,\alpha)\to(\D\otimes A,\gamma\otimes\alpha)
\]
satisfies the conditions stated in \cite[Proposition 3.1]{Sza17}, which are identical to the ones stated in the later article \cite[Proposition 4.3]{Sza21}.
The conclusion of the latter yields that $\1_{\D}\otimes\id_A$ is asymptotically unitarily equivalent to a cocycle conjugacy, which is condition \ref{prop:embeddingD2c:3}.

\ref{prop:embeddingD2c:3}$\Rightarrow$\ref{prop:embeddingD2c:1}:
Let the cocycle conjugacy $(\f,\bbu)$ and the map $w$ be given as in the statement.
Define a point-norm continuous path $(\Phi_t)_{t\in[0,\infty)}$ of unital embeddings $\D\to A$ given by
\begin{equation*}
\Phi_t:=\f^{-1} \circ \Ad(w_t) \circ (\1\otimes\id_{\D}).
\end{equation*}
Firstly, note that we have for all $a\in A$ and $d\in\D$ that
\[
\lim_{t\to\infty} \|[\Phi_t(d),a]\| = \lim_{t\to\infty} \| [\Phi_t(d),\f^{-1}(w_t(a\otimes\1)w_t^*)] \| = \lim_{t\to\infty} \|[\1\otimes d,a\otimes\1]\| = 0.
\]
Secondly, observe that
\[
\begin{array}{ccl}
\alpha_g(\Phi_t(d))  &=& \alpha_g\circ\f^{-1}\circ\Ad(w_t)(\1\otimes d) \\
&=& (\f\circ\alpha_{g^{-1}})^{-1}\circ\Ad(w_t)(\1\otimes d) \\
&=& (\Ad(\bbu_{g^{-1}})\circ(\alpha\otimes\gamma)_{g^{-1}}\circ\f)^{-1}\circ\Ad(w_t)(\1\otimes d) \\
&=& \f^{-1}\circ(\alpha\otimes\gamma)_g\circ\Ad(\bbu_{g^{-1}}^*w_t)(\1\otimes d) \\
&=& \f^{-1}\circ\Ad(\bbu_g(\alpha\otimes\gamma)_g(w_t))(\1\otimes\gamma_g(d))
\end{array}
\]
for any $g\in G$ and $d\in\D$.
Using the cocycle condition above, one may infer that
\begin{align*}
\lim_{t\to\infty}\max_{g\in K}\|w_t-\bbu_g(\alpha\otimes\gamma)_g(w_t)\|=0
\end{align*}
for every compact subset $K\subseteq G$, and thus that the path $(\Phi_t)_{t\in[0,\infty)}$ is asymptotically equivariant, i.e.,
\begin{equation*}
\lim_{t\to\infty}\max_{g\in K}\|\Phi_t(\gamma_g(d)) - \alpha_g(\Phi_t(d))\| = 0
\end{equation*}
for all $d\in\D$ and every compact subset $K\subseteq G$.
It follows from the previous two observations that $(\Phi_t)_{t\in[0,\infty)}$ induces a well-defined equivariant unital embedding
\[
\Phi: (\D,\gamma) \to (A_{\mathfrak{c},\alpha}\cap A',\alpha_{\mathfrak{c}}).
\]
This finishes the proof.
\end{proof}

The following notion is a straightforward generalization of \cite[Definition 1.4]{Sch20}:

\begin{definition} 
Let $A$ be a \Cs-algebra with a continuous action $\alpha: G\curvearrowright A$.
Let $\gamma: G\curvearrowright\D$ be a strongly self-absorbing action.
We say that $\alpha$ is \emph{separably $\gamma$-stable} if for every separable \Cs-subalgebra $C\subseteq A$, there exists a separable $\alpha$-invariant \Cs-algebra $B\subseteq A$ containing $C$ and such that $(\alpha\restriction{B})$ is $\gamma$-stable.
\end{definition}

Note that any $\gamma$-stable action on a separable \Cs-algebra is trivially seen to be separably $\gamma$-stable.
There are many more straightforward observations one could make about this notion in analogy to \cite[Section 1]{FS24}, but we shall omit them, as they are not particularly relevant for the remainder of this article.

We end this section by recalling a well-known result of Kasparov about the existence of quasi-central approximate units that are compatible with a group action.

\begin{lemma}[see {\cite[Lemma 1.4]{Kas88}}] \label{lem:Kasparov}
Let $\beta: G\curvearrowright B$ be an action on a $\sigma$-unital \Cs-algebra.
Then for any separable \Cs-subalgebra $D\subseteq\M(B)$, there exists a countable, increasing approximate unit of positive contractions $h_n\in B$ satisfying $h_n=h_n h_{n+1}$ for all $n\geq 1$ and
\[
\lim_{n\to\infty} \|[h_n,d]\|=0,\quad d\in D,
\]
and
\[
\lim_{n\to\infty} \max_{g\in K} \| h_n-\beta_g(h_n) \| = 0,\quad K\subseteq G \text{ compact}.
\]
\end{lemma}


\section{The dynamical folding property}

It is known from \cite[Proposition 1.4]{PW07} that the corona algebra $\Q(A)$ of any $\sigma$-unital \Cs-algebra $A$ has the folding property (see Definition \ref{def:folding}).\
The goal of this section is to show that if $A$ is equipped with a $G$-action $\alpha$, then $\Q(A)$ equipped with the (algebraic) $G$-action induced by $\alpha$ satisfies a dynamical generalization of the folding property introduced in Definition \ref{def:dynamical folding}.

Let us first recall that an action on $A$ induces an action on its corona algebra.

\begin{remark}
Let $A$ be a \Cs-algebra equipped with an action $\alpha:G\curvearrowright A$.\
Then $\alpha$ induces an algebraic action on the multiplier algebra of $A$, which we denote again by $\alpha$, given by
\[
\alpha_g(x)\cdot a = \alpha_g(x\cdot \alpha_g^{-1}(a)) \quad\text{for all } g\in G, x\in\M(A) \text{ and } a\in A.
\]
Consequently, one obtains a quotient action induced by $\alpha$ on the corona algebra $\Q(A)=\M(A)/A$.\
We denote this algebraic action by $\alpha$ as well.
\end{remark}

\begin{definition}\label{def:folding}
Let $A$ be a \Cs-algebra.
One says that $A$ has the \textit{folding property} if for every separable \Cs-algebra $B$, \Cs-subalgebra $B_0\subseteq B$, and $\ast$-homomorphism $\psi:B\to A_{\mathfrak{c}}$ with $\psi(B_0)\subseteq A$, there exists a $\ast$-homomorphism $\f:B\to A$ such that $\f\restriction_{B_0}=\psi\restriction_{B_0}$, and $\ker\f=\ker\psi$.
\end{definition}

We note that the last condition about the kernels was not explicitly recorded in the context of corona algebras in \cite{PW07}, but it can be observed as an outcome of the proof.
Since we are about to generalize this property further and hence reprove it anyway, we shall not explain this in greater detail here.

One can summarize the folding property by saying that for each diagram as below, where the square commutes, there exists a map $\f$ with same kernel as $\psi$ that makes the upper-left triangle commute.
\[\begin{tikzcd}
	{B_0} && A \\
	{B} && {A_{\mathfrak{c}}}
	\arrow["\psi\restriction_{B_0}", from=1-1, to=1-3]
	\arrow[hook, from=1-1, to=2-1]
	\arrow["\iota_\mathfrak{c}", hook, from=1-3, to=2-3]
	\arrow["\f"{description}, from=2-1, to=1-3]
	\arrow["{\psi}"', from=2-1, to=2-3]
\end{tikzcd}\]

The following is a dynamical generalization of the folding property.

\begin{definition}\label{def:dynamical folding}
Let $A$ be a \Cs-algebra with an algebraic action $\alpha:G\curvearrowright A$.\footnote{Recall that the \Cs-subalgebra of $A$ containing every element on which the action $\alpha$ is continuous is denoted by $A_\alpha$.}
One says that $(A,\alpha)$ (or $\alpha$) has the \textit{dynamical folding property} if the following statement holds.
Let $B$ be a separable \Cs-algebra and $B_0,D \subseteq B$ \Cs-subalgebras such that $D$ is equipped with a continuous action $\delta:G\curvearrowright D$.
If $\psi:B\to A_{\mathfrak{c}}$ is a $\ast$-homomorphism such that
\begin{itemize}[leftmargin=*]
\item $\psi(B_0)\subseteq A$, and
\item $\psi\restriction_{D}$ is $\delta$-to-$\alpha_{\mathfrak{c}}$ equivariant with range in $(A_{\alpha})_{\mathfrak{c}}$,
\end{itemize}
then there exists a $\ast$-homomorphism $\f:B\to A$ such that
\begin{itemize}[leftmargin=*]
\item $\f\restriction_{B_0}=\psi\restriction_{B_0}$,
\item $\f\restriction_{D}$ is $\delta$-to-$\alpha$ equivariant, and 
\item $\ker\f=\ker\psi$.
\end{itemize}
\end{definition}

One can once again summarize the above definition by saying that $\f$ fits in the following diagram with commuting squares and equivariant restriction $\psi\restriction_{D}$, in such a way that the upper-left triangle commutes, $\f\restriction_{D}$ is equivariant, and $\ker\f=\ker\psi$.
\[
\begin{tikzcd}
	{B_0} && A \\
	B && {A_{\mathfrak{c}}} \\
	D && {(A_{\alpha})_{\mathfrak{c}}}
	\arrow["{\psi\restriction_{B_0}}", from=1-1, to=1-3]
	\arrow[hook, from=1-1, to=2-1]
	\arrow["{\iota_\mathfrak{c}}", hook, from=1-3, to=2-3]
	\arrow["\f", from=2-1, to=1-3]
	\arrow["\psi", from=2-1, to=2-3]
	\arrow[hook, from=3-1, to=2-1]
	\arrow["{\psi\restriction_{D}}", from=3-1, to=3-3]
	\arrow[hook, from=3-3, to=2-3]
\end{tikzcd}
\]
We shall now collect a few technical prerequisites to prove the main result of this section.
The following result by Michael \cite[Proposition 7.2]{Mic56} is a generalization of the Bartle--Graves selection theorem \cite[Theorem 4]{BG52}.

\begin{theorem}[{see \cite[Proposition 7.2]{Mic56}}]\label{thm:BartleGravesOld}
Let $X$ and $Y$ be complex Banach spaces, and $q:X\to Y$ a linear and continuous surjection.\
Then, for any $M>1$, there exists a (not necessarily linear) continuous map $\rho:Y \to X$ such that
\begin{enumerate}[label=\textup{(\roman*)},ref=\textup{\roman*},leftmargin=*]
\item $q \circ \rho = \id_{Y}$,
\item $\|\rho(y)\| \leq M \cdot \inf\{\|x\| \mid x\in q^{-1}(y)\}$ for all $y\in Y$, \label{bound}
\item $\rho(\lambda y)=\lambda\rho(y)$ for all $\lambda\in\mathbb{C}$ and $y\in Y$.
\end{enumerate}
\end{theorem}

In order to apply Michael's result in an efficient manner for our purposes, we restrict ourselves to the case where $X$ and $Y$ are \Cs-algebras.
Under this assumption, $q$ necessarily behaves like a quotient map, and therefore \eqref{bound} amounts exactly to the condition $\|\rho(y)\| \leq M \|y\|$ for all $y\in Y$.
We are now ready to present an application of Theorem \ref{thm:BartleGravesOld} that will play a crucial role in the main result of this section.\
We remark that its proof is inspired by that of \cite[Lemma 1.2]{PW07}.

\begin{corollary}\label{cor:BartleGraves}
Let $q:B \to Q$ be a surjective $\ast$-homomorphism of \Cs-algebras.\
Suppose, moreover, that $B_0 \subseteq B$ and $Q_0 \subseteq Q$ are \Cs-subalgebras that satisfy $q(B_0)=Q_0$.
Then, for each $M>1$, there exists a (not necessarily linear) continuous map $\rho:Q \to B$ such that
\begin{enumerate}[label=\textup{(\roman*)},ref=\textup{\roman*},leftmargin=*]
\item $q \circ \rho = \id_{Q}$,\label{first}
\item $\rho(Q_0) \subseteq B_0$, \label{second}
\item $\|\rho(y)\| \leq 5M \cdot \|y\|$ for every $y\in Q$,\label{third}
\item $\rho(\lambda y)=\lambda\rho(y)$ for every $\lambda\in\mathbb{C}$ and for every $y\in Q$.\label{fourth}
\end{enumerate}
\end{corollary}
\begin{proof}
Fix $M>1$ and $1<m_1<M$, and then set $m_2:=M/m_1>1$, which satisfies $m_1 \cdot m_2=M$.

First of all, we apply Theorem \ref{thm:BartleGravesOld} to $q\restriction_{B_0}:B_0\to Q_0$, and obtain a continuous map $\rho_0:Q_0\to B_0$ such that
\begin{itemize}[leftmargin=*]
\item $(q\restriction_{B_0}) \circ \rho_0 = \id_{Q_0}$,
\item $\|\rho_0(y)\| \leq m_1 \cdot \|y\|$ for all $y\in Q_0$,
\item $\rho_0(\lambda y)=\lambda\rho_0(y)$ for all $y\in Q_0$ and $\lambda\in\mathbb{C}$.
\end{itemize}
Below, we argue that one can extend $\rho_0$ to a right-inverse of $q$ with the desired properties.

We proceed by defining, as in \cite[Lemma 1.2]{PW07}, a retraction $q_0:Q \to Q_0$, namely, a continuous function $q_0:Q \to Q_0$ such that $q_0\restriction_{Q_0}=\id_{Q_0}$.\
By applying Theorem \ref{thm:BartleGravesOld} to the (Banach space) quotient map $p:Q \to Q/Q_0$, we may find a continuous map $\kappa_0:Q/Q_0 \to Q$ such that
\begin{itemize}[leftmargin=*]
\item $p \circ \kappa_0 = \id_{Q/Q_0}$,
\item $\|\kappa_0(p(y))\| \leq m_2 \cdot \|y\|$ for all $y\in Q$,
\item $\kappa_0(\lambda p(y))=\lambda\kappa_0(p(y))$ for all $y\in Q$ and $\lambda\in\mathbb{C}$.
\end{itemize}
Then, the continuous map $q_0:Q \to Q_0$ given by
\begin{equation*}
q_0(y)=y-\kappa_0(p(y))
\end{equation*}
for all $y\in Q$, is the desired retraction.\
Note that $\|q_0(y)\| \leq 2m_2 \cdot \|y\|$ for all $y\in Q$.

We apply once more Theorem \ref{thm:BartleGravesOld}, this time to $q:B\to Q$, and obtain a continuous map $\kappa:Q \to B$ such that
\begin{itemize}[leftmargin=*]
\item $q \circ \kappa = \id_{Q}$,
\item $\|\kappa(y)\| \leq m_1 \cdot \|y\|$ for all $y\in Q$,
\item $\kappa(\lambda y)=\lambda\kappa(y)$ for all $y\in Q$ and $\lambda\in\mathbb{C}$.
\end{itemize}
Now, the continuous function $\rho:Q\to B$ defined by
\begin{equation*}
\rho(y)=\kappa(y-q_0(y)) + \rho_0(q_0(y))
\end{equation*}
for all $y\in Q$, is the desired right-inverse of $q$.\
We check that all conditions are satisfied.\
First of all, note that
\begin{equation*}
q (\rho(y)) = q(\kappa(y-q_0(y))) + q(\rho_0(q_0(y))) = y-q_0(y)+q_0(y)=y
\end{equation*}
for all $y\in Q$, which establishes Condition \eqref{first}.
Since we know that $q_0:Q\to Q_0$ is a retraction, then $\rho\restriction_{Q_0}=\rho_0$, and Condition \eqref{second} follows immediately.
In order to check that Condition \eqref{third} holds, note that
\begin{align*}
\|\rho(y)\| &\leq \|\kappa(y-q_0(y))\| + \|\rho_0(q_0(y))\| \\
&\leq m_1 \left( \|y-q_0(y)\| + \|q_0(y)\| \right) \\
&\leq m_1 \left( \|y\| + 2\|q_0(y)\| \right) \\
&\leq m_1 \left( \|y\| + 4 m_2 \|y\| \right) \\
&\leq 5 M \|y\|
\end{align*}
for all $y\in Y$.\
Finally, since all maps involved preserve multiplication by scalars, $\rho$ satisfies condition \eqref{fourth}.
\end{proof}

\begin{definition}\label{def:equicontinuous}
Let $A$ and $B$ be \Cs-algebras, and $(\sigma_t)_{t\in[0,\infty)}$ a path of maps from $A$ to $B$.
Then $(\sigma_t)_{t\in[0,\infty)}$ is said to be \textit{equicontinuous} if for every $a\in A$ and $\e>0$, there exists $\delta>0$ such that $\|\sigma_t(a)-\sigma_t(b)\|\leq\e$ for all $t\in[0,\infty)$ and $b\in A$ with $\|a-b\|\leq\delta$.
\end{definition}

\begin{lemma} \label{lem:strict-sums}
Let $A$ be a $\sigma$-unital \Cs-algebra.\
Let $(f_n)_{n\in\N}$ be a sequence of pairwise commuting positive contractions in $A$ with $\sum_{n\in\N}f_n^2=\1$ strictly and such that $f_kf_\ell=0$ whenever $|k-\ell|\geq 2$.
Let $Y\subseteq\N$ be any subset.
Then:
\begin{enumerate}[label=\textup{(\roman*)},leftmargin=*]
\item \label{lem:strict-sums:1}
For every bounded sequence $(a_n)_n$ in $\M(A)$, the series $\sum_{n\in Y} f_na_nf_n$ converges strictly in $\M(A)$.
The map $\Psi_Y:\ell^{\infty}(\N,\M(A))\to\M(A)$ given by
\begin{equation*}
\Psi_Y((a_n)_n) = \sum_{n\in Y}f_na_nf_n
\end{equation*}
is completely positive and contractive.
Furthermore $\Psi_Y(c_0(\N,\M(A)))\subseteq A$.
\item \label{lem:strict-sums:2}
Suppose that $\sigma\in\Aut(A)$ is an automorphism such that $\sum_{n=1}^\infty \|\sigma(f_n)-f_n\|<\infty$.
Then one has
\begin{equation*}
\sigma(\Psi_Y((a_n)_n)) - \Psi_Y((\sigma(a_n))_n) \; \in A
\end{equation*}
for every $(a_n)_n\in\ell^{\infty}(\N,\M(A))$.
\end{enumerate}
\end{lemma}
\begin{proof}
\ref{lem:strict-sums:1} is a direct consequence of \cite[Lemma 3.1]{MT04}.

For \ref{lem:strict-sums:2}, observe that every $n\geq 1$ that
\begin{align*}
&\sigma(f_na_nf_n) - f_n\sigma(a_n)f_n \\
&=  \sigma(f_na_nf_n) - f_n\sigma(a_nf_n)
 + f_n\sigma(a_nf_n) - f_n\sigma(a_n)f_n \\
&= (\sigma(f_n)-f_n)\sigma(a_n f_n) + f_n\sigma(a_n)(\sigma(f_n)-f_n).
\end{align*}
By assumption, we get
\[
\sum_{n\in Y} \|\sigma(f_na_nf_n) - f_n\sigma(a_n)f_n\| \leq 2\|(a_n)_n\|\cdot\sum_{n=1}^\infty \|\sigma(f_n)-f_n\|<\infty
\]
and we may conclude that
\begin{equation*}
\sigma(\Psi_Y((a_n)_n)) - \Psi_Y((\sigma(a_n))_n) = \sum_{n\in Y}\sigma(f_na_nf_n) - \sum_{n\in Y}f_n\sigma(a_n)f_n \; \in A.
\end{equation*}
\end{proof}

The following fact is well-known and follows from a standard argument if one approximates the square root function on the unit interval with polynomial functions.

\begin{lemma} \label{lem:epsilon-delta}
For every $\e>0$, there exists a constant $\delta>0$ with the following property.
Let $A$ be any \Cs-algebra, $a\in A$ a positive contraction.
\begin{enumerate}[label=\textup{(\roman*)},leftmargin=2em]
\item If $x\in A$ is any contraction with $\|[x,a]\|<\delta$, then $\|[x,\sqrt{a}]\|<\e$.
\item If $\alpha\in\Aut(A)$ is an automorphism with $\|a-\alpha(a)\|<\delta$, then $\|\sqrt{a}-\alpha(\sqrt{a})\|<\e$.
\end{enumerate}
\end{lemma}

The following is the main result of this section and represents a dynamical generalization of \cite[Proposition 1.4]{PW07}.

\begin{theorem}\label{thm:dynamical folding}
Let $\alpha:G\curvearrowright A$ be an action on a $\sigma$-unital \Cs-algebra, and equip $\Q(A)$ with the algebraic action induced by $\alpha$.\
Then $(\Q(A),\alpha)$ has the dynamical folding property.
\end{theorem}
\begin{proof}
Fix a separable \Cs-algebra $B$ with \Cs-subalgebras $B_0,D\subseteq B$ such that $D$ is equipped with a continuous $G$-action $\delta$.
Let $\psi:B\to\Q(A)_{\mathfrak{c}}$ be a $\ast$-homomorphism such that $\psi(B_0)\subset\Q(A)$, and $\psi\restriction_{D}$ is equivariant with range in $(\Q(A)_\alpha)_{\mathfrak{c}}$.

Consider the quotient map $\omega:\M(A)\to\Q(A)$, and note that $\omega(\M(A)_{\alpha})=\Q(A)_{\alpha}$.\
By Corollary \ref{cor:BartleGraves} (applied with $M=2$), there exists a (not necessarily linear) continuous map $\zeta:\Q(A)\to\M(A)$ such that
\begin{itemize}[leftmargin=*]
\item $\omega \circ \zeta = \id_{\Q(A)}$,
\item $\zeta(\Q(A)_{\alpha}) \subseteq \M(A)_{\alpha}$,
\item $\|\zeta(y)\| \leq 10 \|y\|$ for all $y\in \Q(A)$,
\item $\zeta(\lambda y)=\lambda\zeta(y)$ for all $\lambda\in\mathbb{C}$ and $y\in \Q(A)$.
\end{itemize}
The map
\begin{equation*}
\omega_{*}:\C_b([0,\infty),\M(A))\to\C_b([0,\infty),\Q(A)),\quad \omega_{*}(f)(t)=\omega(f(t))
\end{equation*}
is surjective, and one can see that the continuous map given by
\[
\zeta_{*}:\C_b([0,\infty),\Q(A))\to\C_b([0,\infty),\M(A)), \quad \zeta_{*}(h)(t)=\zeta(h(t))
\]
enjoys the following properties,
\begin{itemize}[leftmargin=*]
\item $\omega_{*} \circ \zeta_{*} = \id_{\C_b([0,\infty),\Q(A))}$,
\item $\zeta_{*}(\C_b([0,\infty),\Q(A)_{\alpha})) \subseteq \C_b([0,\infty),\M(A)_{\alpha})$,
\item $\|\zeta_{*}(h)\| \leq 10 \|h\|$ for all $h\in \C_b([0,\infty),\Q(A))$,
\item $\zeta_{*}(\lambda h)=\lambda\zeta_{*}(h)$ for all $\lambda\in\mathbb{C}$ and $h\in \C_b([0,\infty),\Q(A))$.
\end{itemize}
The quotient map $\pi_{\mathfrak{c}}: \C_b([0,\infty),\Q(A)) \to \Q(A)_{\mathfrak{c}}$ is of course surjective.\
Moreover, after identifying $(\Q(A)_{\alpha})_{\mathfrak{c}}$ as a \Cs-subalgebra of $\Q(A)_{\mathfrak{c}}$, we have that $\pi_{\mathfrak{c}}(\C_b([0,\infty),\Q(A)_{\alpha}))=(\Q(A)_{\alpha})_{\mathfrak{c}}$.\
Using Corollary \ref{cor:BartleGraves} again, we get a continuous map $\rho_{\mathfrak{c}}:\Q(A)_{\mathfrak{c}} \to \C_b([0,\infty),\Q(A))$ such that
\begin{itemize}[leftmargin=*]
\item $\pi_{\mathfrak{c}} \circ \rho_{\mathfrak{c}} = \id_{\Q(A)_{\mathfrak{c}}}$,
\item $\rho_{\mathfrak{c}}((\Q(A)_{\alpha})_{\mathfrak{c}}) \subseteq \C_b([0,\infty),\Q(A)_{\alpha})$,
\item $\|\rho_{\mathfrak{c}}(y)\| \leq 10 \|y\|$ for all $y\in \Q(A)_{\mathfrak{c}}$,
\item $\rho_{\mathfrak{c}}(\lambda y)=\lambda\rho_{\mathfrak{c}}(y)$ for all $\lambda\in\mathbb{C}$ and $y\in \Q(A)_{\mathfrak{c}}$.
\end{itemize}
Hence, we may define an equicontinuous bounded map
\begin{equation*}
\sigma: [0,\infty) \times B \to \M(A), \quad (t,b) \mapsto \sigma_t(b):=\zeta\big( \rho_{\mathfrak{c}}(\psi(b))(t) \big),
\end{equation*}
satisfying the following properties,
\begin{enumerate}[label=\textup{(\roman*)},ref=\textup{\roman*},leftmargin=*]
\item $\omega \circ \sigma_t(b) = \rho_{\mathfrak{c}}(\psi(b))(t)$ for all $b\in B$ and $t\in [0,\infty)$, \label{no1}
\item $[t\mapsto\sigma_{t}(d)] \in \C_b([0,\infty),\M(A)_{\alpha})$ for all $d\in D$, \label{no2}
\item $\sigma_{t}(\lambda b)=\lambda\sigma_{t}(b)$ for all $t\in [0,\infty)$, $\lambda\in\mathbb{C}$ and $b\in B$. \label{no3}
\end{enumerate}
From conditions \eqref{no1} and \eqref{no2} one derives a key property of the path of maps $(\sigma_t)_t$.\
We know that $\omega \circ \sigma_t$ is an asymptotic $\ast$-homomorphism from $B$ to $\Q(A)$ in the sense that it becomes approximately linear, $\ast$-preserving, and multiplicative as $t\to\infty$.\
It is also asymptotically equivariant with respect to $\delta$ and $\alpha$ when restricted to $D\subseteq B$.
Since $\psi(B_0)$ is in the constant part of $\Q(A)_{\mathfrak c}$, we may furthermore conclude that
\begin{equation} \label{eq:sigma-on-B0}
\lim_{t\to\infty} \omega\circ\sigma_t(b) = \psi(b),\quad b\in B_0.
\end{equation}
For the rest of the proof, we choose an increasing sequence of compact subsets $K_n\subseteq G$ such that $\bigcup_n K_n=G$ and an increasing sequence of finite sets $F_n\subseteq B$ such that $F:=\bigcup_n F_n$ is dense in $B$, $F\cap B_0$ is dense in $B_0$, and $F \cap D$ is dense in $D$.

We proceed to find an unbounded increasing sequence $(t_n)_n$ in $[0,\infty)$ such that for all $a\in B$, the difference between $\sigma_{t_n}(a)$ and $\sigma_{t_{n+1}}(a)$ converges to zero as $n\to\infty$.\
(Although this resembles \cite[Lemma 3.3]{MT04}, we write out the proof for the reader's convenience.)
For every element $a\in B$, the function $t\mapsto\sigma_t(a)$ is continuous.\
Hence, for every $n\in\N$, one may find a finite increasing sequence $(s_i^{(n)})_{i=0}^{k_n} \subset [n,n+1]$ with $s_0^{(n)}=n$ and $s_{k_n}^{(n)}=n+1$, such that
\begin{equation*}
\max_{a\in F_n} \max_{s,t\in[s_{i}^{(n)},s_{i+1}^{(n)}]} \|\sigma_s(a)-\sigma_t(a)\| \leq 2^{-n}.
\end{equation*}
Set $r_n:=\sum_{j=0}^{n-1} k_j$ for $n\geq1$, and $r_0=0$.\
Define an increasing sequence $(t_i)_{i\in\N}$ by concatenating the finite sequences above, i.e., by setting $t_0:=0$, and $t_{i+r_n}:=s_i^{(n)} \in [n,n+1]$ for $0 < i \leq k_n$ and $n\in\N$.\
By construction, one has that
\begin{equation}\label{repar0}
\lim_{n\to\infty} \max_{a\in F_{n}} \max_{r_n \leq \ell \leq r_{n+1}} \max_{s,t\in[t_{\ell},t_{\ell+1}]} \|\sigma_{s}(a)-\sigma_t(a)\|=0.
\end{equation}
By equicontinuity of $(\sigma_t)_t$ we may conclude that the condition holds for any $a\in B$.\
To simplify notation, we denote $\sigma_{t_\ell}(a)$ by $\sigma_\ell(a)$ for every $\ell\geq 1$.

We consider the increasing sequence of norm-compact subsets of $\M(A)$ given by
\begin{align*}
 C_n:=& \{ \zeta(\psi(b)) \mid b\in F_n\cap B_0\} \\
&\cup \{ \sigma_\ell(a) \mid 0\leq \ell \leq n+1,\; a\in F_n\}\\
&\cup \{ \sigma_\ell(a)-\zeta(\psi(b)) \mid 0\leq \ell \leq n,\; a\in F_{n}\}\\
&\cup\{ \sigma_\ell(a)\sigma_\ell(b)-\sigma_\ell(ab) \mid 0\leq \ell \leq n,\; a,b\in F_n\} \\
&\cup \{ \sigma_\ell(a)+\lambda\sigma_\ell(b)-\sigma_\ell(a+\lambda b) \mid 0\leq \ell \leq n,\; a,b\in F_n,\; \lambda\in\mathbb{C}, |\lambda|\leq n \} \\
& \cup \{ \sigma_\ell(a)^*-\sigma_\ell(a^*) \mid 0\leq \ell \leq n,\; a\in F_n \} \\
& \cup \{ \alpha_{g}(\sigma_\ell(d))-\sigma_\ell(\delta_g(d)) \mid 0\leq \ell \leq n,\; d\in D\cap F_n,\; g\in K_n\}.
\end{align*}
We remark that compactness of the last component in the definition of $C_n$ follows from the fact that $\sigma_\ell(d)$ is an $\alpha$-continuous element of $\M(A)$ by condition \eqref{no2} for all $\ell\geq 1$ and $d\in D$.

By Lemma~\ref{lem:Kasparov} there exists an approximate unit $(e_n)_{n\in\N}$ of $A$ such that $e_0=0$, and $e_{n+1}e_n=e_n$, and
\begin{equation}\label{eq:approx_invariant}
\lim_{n\to\infty} \Big( \max_{g\in K_n} \|\alpha_{g}(e_n)-e_n\| + \max_{x\in C_n} \|[e_n,x]\| \Big) = 0.
\end{equation}
Let $x\in\M(A)$ be any element.
Given that $(e_n)_n$ is an increasing approximate unit, we have
\begin{equation*}
\|\omega(x)\|=\lim_{n\to\infty} \|(\1-e_n)^{1/2} x\| = \lim_{n\to\infty}\lim_{m\to\infty} \|(e_m-e_n)^{1/2} x\|.
\end{equation*}
In particular, after passing to a subsequence of $(e_n)_n$, we may ensure that
\begin{equation}\label{eq:Cn}
\lim_{n\to\infty} \max_{x\in C_n} \left| \|(\1-e_{n-2})^{1/2} x\| - \|\omega(x)\| \right| = 0
\end{equation}
and
\begin{equation} \label{eq:choice-en+1}
\lim_{n\to\infty} \inf_{m\geq n} \max_{x\in C_n} \left| \|\omega(x) \| - \|(e_{m}-e_{n-1})^{1/2} x \| \right| = 0.
\end{equation}
Note that arbitrary subsequences of $(e_n)_n$ still satisfy conditions \eqref{eq:approx_invariant}, \eqref{eq:Cn} and \eqref{eq:choice-en+1}. 
For every $n\geq 1$, we apply Lemma \ref{lem:epsilon-delta} and find a constant $\delta_n$ satisfying the condition stated there for $2^{-n}$ in place of $\e$.
We may assume $\delta_{n+1}\leq\delta_n\leq 2^{-n}$ for all $n$.
After passing to a subsequence of $(e_n)_n$ again,  we may assume for all $n\geq 1$ that
\[
\max_{g\in K_n}\|\alpha_g(e_n)-e_n\|\leq \delta_{n+1}/2
\]
and
\[
\|[e_n,x]\|\leq \frac{\delta_n}{2(1+\|x\|)},\quad x\in C_n.
\]
We shall denote $e_0=e_{-1}=0$ and set $f_n=(e_{n}-e_{n-1})^{1/2}$ for all $n\geq 0$.
Due to the equation $e_ne_{n+1}=e_n$ and the choice of $\delta_n$, we observe the following properties:
\begin{enumerate}[label=\textup{(\arabic*)},ref=\textup{\arabic*},leftmargin=*]
\item $f_mf_n = 0$ if $|m-n|\geq2$, \label{i1}
\item $\max_{\ell\leq n+1}\|[f_n,\sigma_\ell(a)]\|\leq 2^{-n}\cdot 100\|a\|$\footnote{Here the constant 100 comes the fact that the composition of $\zeta$ with $\rho_{\mathfrak c}$ increases the norm of an element in $\Q(A)$ by a factor of at most $10\cdot 10$.} for all $n\geq 1$ and $a\in F_n$, \label{i2}
\item $\|[f_n,\zeta(\psi(b))]\|\leq 2^{-n}\cdot 10\|b\|$ for all $n\geq 1$ and $b\in F_n\cap B_0$, \label{i2-prime}
\item $\max_{g\in K_n}\|\alpha_g(f_n)-f_n\|\leq 2^{-n}$ for all $n\geq 1$, \label{i3}
\item $\sum_{n=1}^\infty f_n^2=\1$ strictly, \label{i4}
\item $f_n ( f_{n-1}^2+f_n^2+f_{n+1}^2 )=f_n$ for all $n\geq 1$. \label{i5}
\end{enumerate}
We shall now apply Lemma \ref{lem:strict-sums} to this choice of the sequence $(f_n)_n$.
Since $(\omega \circ \sigma_t)_t$ is approximately $\ast$-homomorphic on $B$ and equivariant on $D$ as $t\to\infty$, we can conclude with Lemma \ref{lem:strict-sums} and condition \eqref{eq:Cn} (using $f_n=f_n (1-e_{n-2})^{1/2}$) that
\begin{enumerate}[label=\textup{(\arabic*)},ref=\textup{\arabic*},leftmargin=*,resume]
\item $\sum_{n=1}^\infty f_n \left(\sigma_n(a)\sigma_n(b)-\sigma_n(ab)\right) f_n \in A$ for all $a,b\in B$, \label{lim1}
\item $\sum_{n=1}^\infty f_n \left(\sigma_n(a)+\lambda\sigma_n(b)-\sigma_n(a+\lambda b)\right)f_n \in A$ for all $a,b\in B$ and $\lambda\in\mathbb{C}$, \label{lim2}
\item $\sum_{n=1}^\infty f_n \left(\sigma_n(a)^*-\sigma_n(a^*)\right)f_n \in A$ for all $a\in B$, \label{lim3}
\item $\sum_{n=1}^\infty f_n\left(\alpha_{g}(\sigma_n(d))-\sigma_n(\delta_g(d)) \right)f_n \in A$ for all $g\in G$ and $d\in D$. \label{lim4}
\end{enumerate}
Let $\Phi:B\to\M(A)$ be the map given by
\begin{equation*}
\Phi(b) = \sum_{n=1}^\infty f_n\sigma_n(b)f_n,\quad b\in B,
\end{equation*}
and denote by $\f=\omega\circ\Phi:B\to\Q(A)$ its composition with the quotient map onto the corona algebra.
We show that $\f$ is a $\ast$-homomorphism with $\f\restriction_{B_0}=\psi\restriction_{B_0}$ and such that $\f\restriction_{D}$ is equivariant.
We prove this by verifying that $\Phi$ has these properties modulo $A$.
Since $\Phi$ is continuous (as a composition of continuous maps), it suffices to prove these properties on dense subsets of $B$, $B_0$ and $D$, respectively.
We start with multiplicativity.\
We shall use the notation $a \equiv b$ to say that two elements $a,b\in\M(A)$ agree modulo $A$.
We have for every $a,b\in F$
\begin{align*}
\Phi(a)\Phi(b) &= \left(\sum_{n=1}^\infty f_n \sigma_n(a)f_n \right)\left(\sum_{n=1}^\infty f_n \sigma_n(b)f_n \right)\\
&\ontop{\eqref{i1}}{=} \sum_{n=1}^\infty f_n \sigma_n(a)f_n \left(\sum_{j=n-1}^{n+1} f_j \sigma_j(b)f_j \right)  \\
&\ontop{\eqref{i2}}{\equiv} \sum_{n=1}^\infty f_n\big( \sigma_n(a)\sigma_{n-1}(b) f_{n-1}^2 + \sigma_n(a)\sigma_n(b) f_n^2 + \sigma_n(a)\sigma_{n+1}(b) f_{n+1}^2 \big) f_n
\end{align*}
If we consider condition \eqref{repar0} with Lemma \ref{lem:strict-sums}, we conclude
\begin{align*}
\Phi(a)\Phi(b) &\equiv \sum_{n=1}^\infty f_n \sigma_n(a)\sigma_n(b)(f_{n-1}^2+f_n^2+f_{n+1}^2)f_n  \\
&\ontop{\eqref{i5}}{\equiv} \; \sum_{n=1}^\infty f_n  \sigma_n(a)\sigma_n(b) f_n \\
&\ontop{\eqref{lim1}}{\equiv} \sum_{n=1}^\infty f_n \sigma_n(ab) f_n \ = \ \Phi(ab).
\end{align*}
We proceed analogously to show that $\Phi$ is linear and $\ast$-preserving modulo $A$.
For all $a,b\in B$ and $\lambda\in\mathbb C$ we have
\begin{align*}
\Phi(a+ \lambda b)-\Phi(a) - \lambda\Phi(b) = \sum_{n=1}^\infty f_n \left( \sigma_n(a+ \lambda b) - \sigma_n(a)- \lambda\sigma_n(b) \right)f_n \ontop{\eqref{lim2}}{\equiv} 0
\end{align*}
and
\begin{align*}
\Phi(a^*)-\Phi(a)^* = \sum_{n=1}^\infty f_n \left( \sigma_n(a^*) - \sigma_n(a)^*\right)f_n \ontop{\eqref{lim3}}{\equiv} 0.
\end{align*}
Let $b\in B_0\cap F$.
We estimate
\[
\limsup_{n\to\infty} \|(\1-{e_{n-2}})^{1/2}(\sigma_n(b)-\zeta(\psi(b)))\| \ontop{\eqref{eq:Cn}}{\leq} \limsup_{n\to\infty} \|\omega(\sigma_n(b)-\zeta(\psi(b)) )\|  \ontop{\eqref{eq:sigma-on-B0}}{=} 0.
\]
Using Lemma \ref{lem:strict-sums} again, we observe that
\begin{align*}
\Phi(b) &= \sum_{n=1}^\infty f_n \sigma_n(b)f_n \\
&=\sum_{n=1}^\infty f_n \zeta(\psi(b)) f_n + \sum_{n=1}^\infty f_n \big( \sigma_n(b)-\zeta(\psi(b)) \big) f_n \\
&\ontop{\eqref{i5}}{=} \sum_{n=1}^\infty f_n \zeta(\psi(b)) f_n + \sum_{n=1}^\infty f_n (\1-e_{n-2})^{1/2}\big( \sigma_n(b)-\zeta(\psi(b)) \big) f_n \\
&\equiv \sum_{n=1}^\infty f_n \zeta(\psi(b)) f_n \ \ontop{\eqref{i2-prime}}{\equiv} \ \sum_{n=1}^\infty f_n^2\zeta(\psi(b)) \ \ontop{\eqref{i4}}{=} \ \zeta(\psi(b)).
\end{align*}
This proves that $\psi(b)=\f(b)$ for all $b\in F\cap B_0$ and hence also for all $b\in B_0$ by continuity.
Finally, let us show that $\Phi\restriction_D$ is equivariant modulo $A$.
Let $d\in D\cap F$ and $g\in G$.
Using Lemma \ref{lem:strict-sums} and condition \eqref{i3}, we compute
\begin{align*}
\Phi(\delta_g(d)) - \alpha_g(\Phi(d)) &= \sum_{n=1}^\infty f_n\sigma_n(\delta_g(d))f_n - \alpha_g\left( \sum_n f_n\sigma_n(d)f_n \right) \\
&\equiv \sum_{n=1}^\infty f_n\sigma_n(\delta_g(d))f_n -  \sum_{n=1}^\infty f_n \alpha_g(\sigma_n(d)) f_n \\
 &= \sum_{n=1}^\infty f_n\left(\sigma_n(\delta_g(d))-\alpha_g(\sigma_n(d))\right)f_n \ontop{\eqref{lim4}}{\equiv} 0.
\end{align*}
Thus, we have established that $\f$ is a $\ast$-homomorphism that is equivariant when restricted to $D$ and equal to $\psi$ when restricted to $B_0$.

Let us now check that $\ker\f = \ker\psi$.
It follows from the fact that a right inverse obtained from Corollary \ref{cor:BartleGraves} necessarily maps zero to zero that $\ker\psi\subseteq\ker\f$.\
For the opposite inclusion, choose $b\in\ker\f$.\
Firstly, we point out as a direct consequence of the above chain of computations (used to verify that $\f$ is a $*$-homomorphism) that $0=\f(b^*b)$ is represented by the multiplier $\sum_{k=1}^\infty f_k\sigma_k(b)\sigma_k(b)^* f_k$, hence the series defines an element in $A$.
Thus
\[
\|f_n\sigma_n(b)\sigma_n(b)^* f_n\| \leq \|(\1-e_{n-2}) \Big(\sum_{k=1}^\infty f_k\sigma_k(b)\sigma_k(b)^*f_k\Big) (\1-e_{n-2})\| \xrightarrow{n\to\infty}0.
\]
Thus $\|f_n\sigma_n(b)\|\xrightarrow{n\to\infty}0$.\
If we apply \eqref{eq:choice-en+1}, we get for every $a\in F$ that
\begin{align*}
\limsup_{n\to\infty}\|f_n\sigma_n(a)\|=\limsup_{n\to\infty} \|\omega(\sigma_n(a))\|.
\end{align*}
By equicontinuity and the fact that $F$ is dense in $B$, this equation persists for all $a\in B$.
We may therefore conclude $\|\omega(\sigma_n(b))\|\xrightarrow{n\to\infty}0$.\
Given condition \eqref{repar0}, this implies $\|\omega(\sigma_t(b))\|\xrightarrow{t\to\infty}0$ or equivalently $\rho_{\mathfrak{c}}(\psi(b))\in\C_0([0,\infty),\Q(A))$.\
Therefore, we have that $\psi(b)=\pi_{\mathfrak{c}}(\rho_{\mathfrak{c}}(\psi(b)))=0$, and hence that $\ker\f=\ker\psi$.
\end{proof}

\begin{notation}
Let $\alpha: G\curvearrowright A$ be an algebraic action on a \Cs-algebra.
Given a \Cs-subalgebra $C\subseteq A$, we write
\[
A\cap_\alpha C' := A\cap\Big(\bigcup_{g\in G} \alpha_g(C) \Big)'.
\]
This is clearly the largest $\alpha$-invariant \Cs-subalgebra of $A\cap C'$.
\end{notation}

Before we move on to the next section, we make a technical observation that is needed therein, which we obtain as a consequence of the dynamical folding property.

\begin{lemma} \label{lem:asymptotic-to-genuine-equivalence}
Let $\alpha: G\curvearrowright A$ be an algebraic action with the dynamical folding property.
Let $C\subseteq A$ be a separable \Cs-subalgebra and $\delta: G\curvearrowright D$ a continuous action on a separable \Cs-algebra.
Suppose that $\f ,\psi: (D,\delta)\to (A\cap_\alpha C',\alpha)$ are two equivariant $*$-homomorphisms that are asymptotically $G$-unitarily equivalent, i.e., there exists a norm-continuous path $w: [0,\infty)\to\mathcal{U}\big( \1+(A\cap_\alpha C')_\alpha \big)$ with
\[
\psi(d)=\lim_{t\to\infty} w_t\f(d)w_t^*,\quad \lim_{t\to\infty}\max_{g\in K} \|w_t-\alpha_g(w_t)\|=0
\]
for all $d\in D$ and every compact set $K\subseteq G$.
Then $\f$ and $\psi$ are $G$-unitarily equivalent, i.e., there exists a unitary $u\in\mathcal{U}(\1+(A\cap_\alpha C')^\alpha)$ such that $\psi=\Ad(u)\circ\f$.
\end{lemma}
\begin{proof}
The given unitary path induces a unitary element
\[
\bar{w}\in \1+\big( (A\cap_\alpha C')_\alpha \big)_\mathfrak{c}^{\alpha_\mathfrak{c}} \subseteq \1+\big( \big( (A_{\alpha})_\mathfrak{c} \big)^{\alpha_\mathfrak{c}} \cap C' \big)
\]
such that $\psi=\Ad(\bar{w})\circ\f$ when viewed as equivariant $*$-homomorphisms from $D$ to $(A_\alpha)_\mathfrak{c}\cap C'$.
We shall apply the dynamical folding property with the choice $B_0=\Cs(\f(D)\cup\psi(D)\cup C)$, the algebra $\Cs(\bar{w}-\1)$ in place of $D$ equipped with the trivial $G$-action, and $B$ the \Cs-algebra generated by these two sets.
This yields a $\ast$-homomorphism $\kappa: B\to A$ satisfying $\kappa\circ\f=\f$, $\kappa\circ\psi=\psi$ and $\kappa(c)=c$ for all $c\in C$.
The element $u=\1+\kappa(\bar{w}-\1)$ then defines a unitary in $\1+A^\alpha$ because $\kappa$ is equivariant on $\Cs(\bar{w}-\1)$.
Moreover $u$ commutes with $C$ because $\kappa$ is the identity map on $C$.
Next, we have that 
\[
\Ad(u)\circ\f = \kappa\circ\big(\Ad(\bar{w})\circ\f)=\kappa\circ\psi=\psi.
\]
Since $A^\alpha\cap C'=(A\cap_\alpha C')^\alpha$, the element $u$ satisfies the desired property.
\end{proof}


\section{Absorption properties}

In this section we prove the main result of the article, Theorem \ref{thm:saturation}, which generalizes \cite[Theorem 2.5]{FS24} to the dynamical setting.
The following generalizes $D$-saturation, which we have recalled in the introduction, to the dynamical setting.

\begin{definition}\label{def:gamma-saturated}
Let $\alpha:G\curvearrowright A$ be an algebraic action on a unital \Cs-algebra, and $\gamma:G\curvearrowright D$ an action on a separable unital \Cs-algebra.\
We say that $(A,\alpha)$ is \textit{$\gamma$-saturated} if for every separable  \Cs-subalgebra $C\subseteq A$, there exists an equivariant unital embedding $(D,\gamma)\to(A\cap_\alpha C',\alpha)$.
\end{definition}

The next proposition represents one of the main technical steps towards the proof of the main result.
In a special case, it was originally obtained in the first-named author's master’s thesis \cite{Li-thesis}.

\begin{proposition}\label{prop:D2corona}
Let $\alpha:G\curvearrowright A$ be an action on a $\sigma$-unital, non-unital \Cs-algebra and $\gamma:G\curvearrowright\D$ a strongly self-absorbing action.\
Suppose that $\alpha$ is separably $\gamma$-stable.
Then for every separable \Cs-subalgebra $C\subseteq\Q(A)$, there exists an equivariant unital $\ast$-homomorphism
\begin{equation*}
(\D^{(2)},\gamma^{(2)}) \to (\Q(A)\cap_\alpha C',\alpha).
\end{equation*}
\end{proposition}
\begin{proof}
Fix a separable \Cs-subalgebra $C\subseteq\Q(A)$.
Due to \cite[Lemma 2.4]{FS24} and the fact that $\alpha$ is separably $\gamma$-stable, there exists a separable nondegenerate $\alpha$-invariant \Cs-subalgebra $A_1\subseteq A$ such that $(\alpha\restriction{A_1})$ is $\gamma$-stable and under the canonical embedding $\Q(A_1)\subseteq\Q(A)$, we have $C\subseteq\Q(A_1)$.
Considering what the claim says, it suffices to prove it for $(A_1,\alpha)$ in place of $(A,\alpha)$.
In other words, we may assume without loss of generality that $A$ is separable and that $\alpha$ is $\gamma$-stable.

In order to find a $\ast$-homomorphism as in the claim, we shall appeal to the universal property of $(\D^{(2)},\gamma^{(2)})$ as specified in Remark \ref{rem:universal property}.\
In other words, we want to find two equivariant c.p.c.\ order zero maps $\mu_0,\mu_1:(\D,\gamma)\to(\Q(A)\cap_\alpha C',\alpha)$ with commuting ranges such that $\mu_0(\1)+\mu_1(\1)=\1$.

Pick a sequence $(c_n)_{n\geq 1}$ in the unit ball of $\M(A)$ whose image under the quotient map $\omega:\M(A)\to\Q(A)$ is dense in the unit ball of $C\subseteq\Q(A)$.\
Also fix an increasing sequence of finite subsets $F_n$ of the unit ball of $\D$ with dense union $F:=\bigcup_n F_n$, and an increasing sequence of compact subsets $K_n\subseteq G$ whose union is $G$.
For every $n\geq 1$, we apply Lemma \ref{lem:epsilon-delta} and find a constant $\delta_n>0$ such that for any positive contraction $a\in A$ one has the following two implications:
\begin{enumerate}[label={$\bullet$},leftmargin=*]
\item If $\|[a,x]\|<\delta_n$ for some contraction $x\in A$, then $\|[\sqrt{a},x]\|<2^{-n}$.
\item If $\|\alpha(a)-a\|<\delta_n$ for some automorphism $\alpha\in\Aut(A)$, then $\|\alpha(\sqrt{a})-\sqrt{a}\|<2^{-n}$.
\end{enumerate}
We may assume that $\delta_{n+1}\leq\delta_n\leq 2^{-n}$ for all $n$.
By Lemma~\ref{lem:Kasparov} there exists an increasing approximate unit $(e_n)_{n\geq 1}$ of $A$ such that  $e_{n+1}e_n=e_n$ for all $n\geq 1$, and
\[
\max_{\ell\leq {n+1}} \|[c_\ell, e_n]\| + \max_{g\in K_{n+1}} \|\alpha_g(e_n)-e_n\| < \frac{\delta_{n+1}}{2} \quad\text{for all } n\geq 1.
\]
Let us define $e_{-1}=e_0:=0$ and $f_n=(e_{n}-e_{n-1})^{1/2}$ for all $n\geq 0$.
With the triangle inequality, it follows that 
\begin{align*}
\max_{\ell\leq n}   \|[c_\ell, e_n-e_{n-1}]\|  &+ \max_{g\in K_n} \|\alpha_g(e_n-e_{n-1})-(e_n-e_{n-1})\| \\
&< \max_{\ell\leq n} \|[c_\ell, e_n]\|+  \max_{g\in K_n} \|\alpha_g(e_n)-e_n\|  + \frac{\delta_n}{2} \\
&\leq \max_{\ell\leq n+1} \|[c_\ell, e_n]\|+  \max_{g\in K_{n+1}} \|\alpha_g(e_n)-e_n\|  + \frac{\delta_n}{2} \\
&<\frac{\delta_{n+1}}{2} + \frac{\delta_n}{2} \leq \delta_{n}
\end{align*}
for all $n\geq 0$.
Note that $(e_n-e_{n-1})$ is a positive contraction.
Due to our choice of the constants $\delta_n$, this condition implies that 
\begin{enumerate}[label=\textup{(\alph*)},ref=\textup{\alph*},leftmargin=*]
\item $\|[f_n,c_\ell]\|< 2^{-n}$ for all $\ell\leq n$ and $n\geq 1$, \label{a}
\item $\max_{g\in K_n}\|\alpha_g(f_n)-f_n\|< 2^{-n}$ for all $n\geq 1$. \label{b}
\end{enumerate}
Moreover, we have by construction that
\begin{enumerate}[resume,label=\textup{(\alph*)},ref=\textup{\alph*},leftmargin=*]
\item $f_mf_n=0$ if $|m-n|\geq 2$, \label{item1}
\item $\sum_{n=1}^\infty f_n^2=\1$ strictly. \label{item3}
\end{enumerate}
Due to Theorem \ref{thm:ssa-action} we may find a sequence of unital $\ast$-homomorphisms $\f_n:\D\to\M(A)$ such that
\begin{enumerate}[resume,label=\textup{(\alph*)},ref=\textup{\alph*},leftmargin=*]
\item $\max_{d\in F_n}\|[\f_n(d),f_m]\|\leq 2^{-n}$ for all $m\leq n+1$ and $n\geq 1$, \label{c}
\item $\max_{d\in F_n}\|[\f_n(d),f_nc_\ell]\|\leq 2^{-n}$ for all $\ell\leq n$ and $n\geq 1$, \label{d}
\item $\max_{d, d'\in F_n}\|[\f_n(d),f_\ell\f_\ell(d')f_\ell]\|\leq 2^{-n}$ for all $\ell< n$ and $n\geq 1$, \label{e}
\item $\max_{d\in F_n} \max_{g\in K_n} \|f_n (\alpha_g(\f_n(d)) - \f_n(\gamma_g(d))) f_n\| \leq 2^{-n}$ for all $n\geq 1$. \label{f}
\end{enumerate}
Define two maps $\Psi_0,\Psi_1:\D\to\M(A)$ given by
\begin{equation*}
\Psi_i(d) = \sum_{n=0}^\infty f_{2n+i} \f_{2n+i}(d) f_{2n+i},\quad d\in\D,\ i=0,1.
\end{equation*}
By Lemma \ref{lem:strict-sums}, we know that $\Psi_0$ and $\Psi_1$ are c.p.c.\ maps.
The formula $\Psi_0(\1)+\Psi_1(\1)=\1$ evidently holds by construction due to condition \eqref{item3}.
As before, we shall use the notation $a \equiv b$ to say that two elements $a,b\in\M(A)$ agree modulo $A$.
For every $g\in G$ and $d\in\D$, Lemma \ref{lem:strict-sums} furthermore implies with condition \eqref{b} that
\begin{equation*}
\alpha_g(\Psi_i(d)) \equiv \sum_{n=0}^\infty f_{2n+i} \alpha_g(\f_{2n+i}(d)) f_{2n+i},\quad i=0,1.
\end{equation*}
We may therefore infer with condition \eqref{f} that
\begin{equation*}
\alpha_g(\Psi_i(d)) \equiv \Psi_i(\gamma_g(d)) \quad\text{for all } g\in G,\ d\in F \text{ and } i=0,1.
\end{equation*}
Since $F$ is dense in the unit ball of $\D$, we conclude that $\alpha_g(\Psi_i(d)) \equiv \Psi_i(\gamma_g(d))$ holds for all $g\in G$ and $d\in\D$.
In other words, the c.p.c.\ maps given by
\begin{equation*}
\mu_i=\omega \circ \Psi_i: \D \to \Q(A),\quad i=0,1,
\end{equation*}
are equivariant with respect to $\gamma$ and $\alpha$.
To show that the maps $\mu_i$ have range in $\Q(A)\cap C'$, pick $d\in F_n$ for some $n\geq 1$.
Note that for all $j\leq n$,
\begin{equation*}
f_n \f_n(d) f_n c_j \ontop{\eqref{d}}{=}_{2^{-n}} f_n^2 c_j\f_n(d) \ontop{\eqref{a}}{=}_{2^{1-n}} c_j f_n^2\f_n(d) \ontop{\eqref{c}}{=}_{2^{-n}} c_j f_n\f_n(d)f_n.
\end{equation*}
Thus, $\|[f_n \f_n(d) f_n, c_j]\|\leq 2^{2-n}$ for $j\leq n$ and $d\in F_n$, which implies that $[\Psi_i(d),c_j] \in A$ for all $j\geq 1$ and $d\in F$, where $i=0,1$.
Since $F$ and $(\omega(c_n))_n$ are dense in the unit ball of $\D$ and $C$, respectively, it follows for every $d\in \D$ and $c\in C$ that $[\Psi_i(d),c]\in A$, and therefore $\mu_i(\D)\subseteq \Q(A)\cap C'$ for $i=0,1$.
Since we have already shown that $\mu_0$ and $\mu_1$ are equivariant, this yields $\mu_i(\D)\subseteq \Q(A)\cap_\alpha C'$ for $i=0,1$.

In order to show that the ranges of $\mu_0$ and $\mu_1$ commute, let us choose $d, d'\in F_n$.
We note that
\begin{align*}
f_{n+1}\f_{n+1}(d')f_{n+1} \cdot f_n\f_n(d)f_n &\ontop{\eqref{c}}{=}_{2^{-(n+1)}} f_{n+1}^2\f_{n+1}(d')f_n\f_n(d)f_n \\
&\ontop{\eqref{e}}{=}_{2^{-(n+1)}} f_{n+1}^2 f_n\f_n(d)f_n \f_{n+1}(d') \\
&\ontop{\eqref{c}}{=}_{2^{1-n}} f_n\f_n(d)f_n f_{n+1}^2 \f_{n+1}(d') \\
&\ontop{\eqref{c}}{=}_{2^{-(n+1)}} f_n\f_n(d)f_n f_{n+1}\f_{n+1}(d')f_{n+1}.
\end{align*}
In particular, this implies that
\begin{equation*}
\| [f_{n+1}\f_{n+1}(d')f_{n+1} , f_n\f_n(d)f_n] \| \leq 2^{2-n}.
\end{equation*}
By using condition \eqref{item1}, we then have that
\begin{align*}
[\Psi_0(d),\Psi_1(d')] =& \left[ \sum_{n=0}^\infty f_{2n}\f_{2n}(d)f_{2n} \ , \ \sum_{n=0}^\infty f_{2n+1}\f_{2n+1}(d')f_{2n+1} \right] \\
=& \sum_{n=0}^\infty [f_{2n}\f_{2n}(d)f_{2n} , f_{2n+1}\f_{2n+1}(d')f_{2n+1} ] \\
&+ \sum_{n=1}^\infty [f_{2n}\f_{2n}(d)f_{2n} , f_{2n-1}\f_{2n-1}(d')f_{2n-1} ].
\end{align*}
The estimate above implies that these define norm-convergent series in $A$ and thus $[\Psi_0(d),\Psi_1(d')] \in A$ for every $d,d'\in F$.\
Since $F$ was dense in the unit ball of $\D$, the same holds for $d,d'\in\D$, and hence $\mu_0$ and $\mu_1$ have commuting ranges.

To see that $\mu_0$ and $\mu_1$ are order zero maps, we use once again condition \eqref{item1} to see that, for $i=0,1$ and for all $d,d'\in\D$,
\begin{align*}
\Psi_i(d)\Psi_i(d') = \sum_{n=0}^\infty f_{2n+i}\f_{2n+i}(d)f_{2n+i}^2\f_{2n+i}(d')f_{2n+i}.
\end{align*}
If $d,d'\in F_n$ for some $n\geq 1$, then
\begin{align*}
f_{2n+i}\f_{2n+i}(d)f_{2n+i}^2\f_{2n+i}(d')f_{2n+i} &\ontop{\eqref{c}}{=}_{2^{1-(2n+i)}} f_{2n+i}\f_{2n+i}(dd')f_{2n+i}^3.
\end{align*}
We may conclude that
\begin{equation*}
\Psi_i(d)\Psi_i(d') \equiv \Psi_i(ab)\Psi_i(\1)
\end{equation*}
for all $d,d'\in F$, and hence also for all $d,d'\in\D$ .\
To summarize, $\mu_0$ and $\mu_1$ are equivariant c.p.c.\ order zero maps with commuting ranges such that $\mu_0(\1)+\mu_1(\1)=\1$, and thus the universal property recalled in Remark \ref{rem:universal property} yields an equivariant unital $\ast$-homomorphism $(\D^{(2)},\gamma^{(2)}) \to \Q(A)\cap_\alpha C'$.
\end{proof}

We are now ready to prove the main theorem of this article.

\begin{theorem} \label{thm:saturation}
Let $\alpha:G\curvearrowright A$ be an action on a $\sigma$-unital, non-unital \Cs-algebra, and $\gamma:G\curvearrowright\D$ a strongly self-absorbing, unitarily regular action.\
If $\alpha$ is separably $\gamma$-stable, then $(\Q(A),\alpha)$ is $\gamma$-saturated.
\end{theorem}
\begin{proof}
Let $C\subseteq\Q(A)$ be a separable \Cs-subalgebra.
By Proposition \ref{prop:D2corona}, there exists an equivariant unital $\ast$-homomorphism
\begin{equation*}
\kappa: (\D^{(2)},\gamma^{(2)}) \to (\Q(A)\cap_\alpha C', \alpha).
\end{equation*}
It extends to an equivariant unital $\ast$-homomorphism
\begin{equation*}
\kappa_\mathfrak{c}: \D^{(2)}_{\mathfrak{c}} \to (\Q(A)_\alpha)_\mathfrak{c} \quad\text{with}\quad [\kappa_\mathfrak{c}(\D^{(2)}),C]=0
\end{equation*}
that sends a representative function $f:[0,\infty)\to\D^{(2)}$ of an element in $\D^{(2)}_{\mathfrak{c}}$ to $\kappa \circ f \in \C_b([0,\infty),\Q(A)_\alpha\cap C')$.\
This is well defined because whenever $f(t)\xrightarrow{t\to\infty}0$ one gets $\kappa(f(t))\xrightarrow{t\to\infty}0$.
Since the strongly self-absorbing action $\gamma$ is assumed to be unitarily regular, $\gamma^{(2)}$ is $\gamma$-absorbing by Theorem \ref{thm:unitarily-regular}.\
By Proposition \ref{prop:embeddingD2c} there exists an equivariant unital embedding
\begin{equation*}
\theta: (\D,\gamma) \hookrightarrow (\D^{(2)}_{\mathfrak{c},\gamma^{(2)}},\gamma^{(2)}_{\mathfrak{c}}).
\end{equation*}
The resulting composition $\psi=\kappa_\mathfrak{c}\circ\theta$ is a unital equivariant $*$-homomorphism from $(\D,\gamma)$ to $(\Q(A)_\alpha)_\mathfrak{c}$ with $[\psi(\D),C]=0$.
We want to appeal to the dynamical folding property via Theorem~\ref{thm:dynamical folding}.
We make the choice $D=\psi(\D)$ with action induced from $\gamma$, we choose $B_0=\Cs(C,\1)$, and $B=\Cs(B_0\cup D)\subseteq \Q(A)_\mathfrak{c}$.
This allows us to find a $\ast$-homomorphism $\f: B\to\Q(A)$ such that $\f(b)=b$ for all $b\in B_0$ and $\f\circ\psi$ defines an equivariant $*$-homomorphism from $(\D,\gamma)$ to $(\Q(A),\alpha)$.
Since $\psi(\1_{\D})=\1_{\Q(A)}\in B_0$, it follows that $\f\circ\psi$ is unital.
Since $C\subseteq B_0$ and $[\psi(\D),C]=0$, it follows that also $[(\f\circ\psi)(\D),C]=0$.
In other words, we have found an equivariant unital $*$-homomorphism from $(\D,\gamma)$ to $(\Q(A)\cap_\alpha C',\alpha)$.
\end{proof}

We would like to end this section by pointing out that the main result also leads to a uniqueness theorem that generalizes the one observed in \cite[Theorem E]{Farah23}.
We shall deduce this as a consequence of a more general formal observation.

\begin{theorem} \label{thm:uniqueness-gamma-folding}
Let $A$ be a unital \Cs-algebra and $\alpha: G\curvearrowright A$ an algebraic action with the dynamical folding property.
Suppose that $\gamma:G\curvearrowright\D$ is a strongly self-absorbing, unitarily regular action such that $\alpha$ is $\gamma$-saturated.
Let $C\subseteq A$ be a separable \Cs-subalgebra.
Then all equivariant unital $*$-homomorphisms from $(\D,\gamma)$ to $(A\cap_\alpha C',\alpha)$ are mutually $G$-unitarily equivalent.
\end{theorem}
\begin{proof}
Fix $C\subseteq A$ as in the statement.
We first note that due to the assumption, there must exist an equivariant unital $*$-homomorphism from $(\D,\gamma)$ to $(A\cap_\alpha C',\alpha)$ to begin with.
Let $\psi_1$ and $\psi_2$ be two arbitrary ones.
Consider $E=\Cs(\psi_1(\D)\cup\psi_2(\D))$, which is a separable $\alpha$-invariant \Cs-subalgebra of $A$.
Using that $\alpha$ is $\gamma$-saturated, we find a third equivariant unital $*$-homomorphism $\psi_3: (\D,\gamma)\to (A\cap_\alpha (C\cup E)',\alpha)$.
By the universal property of the tensor product, we obtain two equivariant unital $*$-homomorphisms
\[
\theta_1, \theta_2: (\D\otimes\D,\gamma\otimes\gamma)\to (A\cap_\alpha C',\alpha)
\]
defined by $\theta_1(a\otimes b)=\psi_1(a)\psi_3(b)$ and $\theta_2(a\otimes b)=\psi_2(a)\psi_3(b)$ for all $a,b\in\D$.
In other words, we have
\[
\psi_1=\theta_1\circ(\id_\D\otimes\1),\ \psi_2=\theta_2\circ(\id_\D\otimes\1),\ \psi_3=\theta_1\circ(\1\otimes\id_D)=\theta_2\circ(\1\otimes\id_D).
\]
As $\gamma$ is strongly self-absorbing and unitarily regular, it follows from \cite[Theorem 3.15]{Sza18} that the two equivariant unital embeddings
\[
\id_D\otimes\1,\ \1\otimes\id_\D: (\D,\gamma) \to (\D\otimes\D,\gamma\otimes\gamma)
\]
are asymptotically $G$-unitarily equivalent.
This immediately implies with the above equations (via transitivity) that $\psi_1$ and $\psi_2$ are asymptotically $G$-unitarily equivalent as maps into $A\cap_\alpha C'$.
The claim thus follows from Lemma~\ref{lem:asymptotic-to-genuine-equivalence}.
\end{proof}

\begin{corollary}\label{cor:G-unitarily-equivalence}
Let $A$ be a $\sigma$-unital, non-unital \Cs-algebra and $\alpha: G\curvearrowright A$ an action.
We consider the induced algebraic action $\alpha: G\curvearrowright\Q(A)$.
Suppose that $\gamma:G\curvearrowright\D$ is a strongly self-absorbing, unitarily regular action such that $\alpha$ is separably $\gamma$-stable.
Let $C\subseteq\Q(A)$ be any separable \Cs-subalgebra.
Then all equivariant unital $*$-homomorphisms from $(\D,\gamma)$ to $(\Q(A)\cap_\alpha C',\alpha)$ are mutually $G$-unitarily equivalent.
\end{corollary}
\begin{proof}
Combine Theorems~\ref{thm:dynamical folding}, \ref{thm:saturation} and \ref{thm:uniqueness-gamma-folding}.
\end{proof}


\section{Concluding remarks}

We would like to point out that our main result has a partial converse that follows from the results in \cite{FS24}, which we shall summarize and recall below.

\begin{theorem} \label{thm:converse}
Suppose that $A$ is a nonzero \Cs-algebra with an action $\alpha: G\curvearrowright A$.
We equip $\K\otimes A$ with the action $\alpha^s:=\id_{\K}\otimes\alpha$.
Denote by $\pi: \M(\K\otimes A)\to\Q(\K\otimes A)$ the (equivariant) quotient map.
Then there exist
	\begin{enumerate}[label=\textup{(\arabic*)}, leftmargin=*]
	\item a unital equivariant $*$-embedding $\Theta\colon \Q(c_0(A))\to \Q(\K\otimes A)$ with $\Theta(\M(A))=\pi(1\otimes\M(A))$, where $\M(A)$ is identified with its isomorphic image inside $\Q(c_0(A))=\ell^{\infty}(\M(A))/c_0(A)$ under the diagonal embedding,
	\item an equivariant u.c.p.\ map $\Psi: \Q(\K\otimes A)\to \Q(c_0(A))$ with $\Psi\circ\Theta=\operatorname{id}$,
	\item a separable unital \Cs-subalgebra $C$ of $\pi(\M(\K)\otimes 1)\cap\Theta(\Q(c_0(A)))'$ (and hence of $\Q(\K\otimes A)^{\alpha^s}$),
	\end{enumerate}
    such that after composing with the quotient map $\Q(c_0(A))\to\M(A)_\infty$, the restriction
	\[
	\bar{\Psi}: \Q(\K\otimes A)\cap C'\to\M(A)_\infty
	\]
	is an equivariant $*$-homomorphism.
\end{theorem}
\begin{proof}
This follows directly from \cite[Theorem 4.4]{FS24} and the construction of the involved maps in `Definition 4.2' therein, which automatically makes all these maps equivariant.
\end{proof}

\begin{corollary} \label{cor:converse}
Let $\alpha: G\curvearrowright A$ be an action on a separable \Cs-algebra.
Let $\gamma: G\curvearrowright\D$ be a strongly self-absorbing action.
If the algebraic action induced by $\alpha^s:=\id_{\K}\otimes\alpha$ on $\Q(\K\otimes A)$ is $\gamma$-saturated, then $\alpha$ is $\gamma$-stable.
\end{corollary}
\begin{proof}
Let $C\subset\Q(\K\otimes A)^{\alpha^s}$ be a separable \Cs-subalgebra as in the conclusion of Theorem \ref{thm:converse}.
Consider the resulting equivariant $*$-homomorphism
\[
\bar{\Psi}: \Q(\K\otimes A)\cap C'\to\M(A)_\infty.
\]
Notice that by construction, one has $\bar{\Psi}(\1_{\K}\otimes a)=a$ for all $a\in\M(A)$.
Since we assumed that $\Q(\K\otimes A)$ with the induced action is $\gamma$-saturated, it follows that there exists a unital equivariant $*$-homomorphism $\D\to\Q(\K\otimes A)\cap (C\cup (\1\otimes A))'$.
The composition with $\bar{\Psi}$ yields a unital equivariant $*$-homomorphism $\D\to\M(A)_\infty\cap A'$.
If we compose this with the quotient map onto $(\M(A)_\infty\cap A')/(\M(A)_\infty\cap A^\perp)$, then we see that the criterion in Theorem \ref{thm:ssa-action} is fulfilled.
Hence $\alpha$ is $\gamma$-absorbing.
\end{proof}

\begin{remark}
If one were to work out the technical proofs to suitably generalize \cite[Proposition 1.7, Corollary 1.10]{FS24} to the dynamical context, then Corollary~\ref{cor:converse} would also generalize in a straightforward way.
Namely we would not have to assume that $A$ is separable and the conclusion would yield that $\alpha$ is separably $\gamma$-stable.
Since we are currently not aware of any interesting applications of this more general observation, we shall not work out its proof here, but nevertheless wish to point out how to get there if it should be of interest to the reader in their own work.
\end{remark}


\bibliography{references.bib}
\bibliographystyle{mybstnum.bst}

\end{document}